\documentclass[10pt,letterpaper]{article}
\usepackage{blindtext,titlefoot}
\usepackage{authblk}
\usepackage{marvosym}
\usepackage[utf8]{inputenc}
\usepackage[english]{babel}
\usepackage{amsmath}
\usepackage{amsfonts}
\usepackage{amssymb}
\usepackage{addfont}
\usepackage{makeidx}
\usepackage{graphicx}
\usepackage{lmodern}
\usepackage{kpfonts}
\usepackage{dsfont}
\usepackage{cancel}
\usepackage{amsthm} 
\usepackage[left=2cm,right=2cm,top=2cm,bottom=2cm]{geometry}
\usepackage{subcaption}
\usepackage{multirow}
\usepackage{float}
\usepackage{url}
\usepackage{hyperref}
\usepackage{multicol}
\usepackage{mathtools, cuted}
\usepackage{lipsum}
\usepackage{booktabs}
\usepackage{comment}
\usepackage{cite}

\providecommand{\nset}[1]{
\mathbb{#1}
}
\providecommand{\ds}[1]{
\displaystyle #1
}

\newcommand{\set}[1]{
\left\{#1\right\}
}

\newcommand{\salg}[1]{
\mathfrak{ #1 }
}

\newcommand{\norm}[1]{
\left\lVert #1 \right\rVert
}
\newcommand{\abs}[1]{
\left\lvert #1 \right\rvert
}

\newcommand{\scp}[1]{
\left\langle #1  \right\rangle
}
\newcommand{\cexp}[2]{
\mathcal{E}( #1 | \salg{#2})
}

\newcommand{\cexpperp}[1]{
\mathcal{E}_{{#1} }^{\perp}
}
\newtheorem{theorem}{ Theorem}[section]
\newtheorem{definition}[theorem]{Definition}
\newtheorem{proposition}[theorem]{Proposition}

\newtheorem{lemma}[theorem]{Lemma}

\usepackage{enumitem} 
\setlist[itemize]{noitemsep} 

\usepackage{titlesec} 
\titleformat{\section}[block]{\large\bfseries\scshape\centering}{\thesection.}{1em}{} 
\titleformat{\subsection}[block]{\large\bfseries\scshape\centering}{\thesubsection.}{1em}{}
\titleformat{\subsubsection}[block]{\large\bfseries\scshape\centering}{\thesubsubsection.}{1em}{} 

\title{\huge\bfseries Almost everywhere continuity of conditional expectations }

\author[,a]{Alberto Alonso \footnote{Email address: alonsoalber@gmail.com}}
\affil[a]{\normalsize Department of Mathematics, Faculty of Science - UNAM, Mexico}

\author[,a]{Fernando Brambila-Paz \footnote{Email address: fernandobrambila@gmail.com}}

\date{}

\begin{document}
\maketitle


\begin{abstract}

A necessary and sufficient condition on a sequence $\set{\mathfrak{A}_n}_{n\in \nset{N}}$ of $\sigma$-subalgebras that assures convergence almost every where of conditional expectations is given.

\textbf{Keywords:} 
\end{abstract}

\section{Introduction}

Let $\salg{A}$ be a $\sigma$-algebra of a set $\nset{X}$ with a probability measure $\mu$. Given a sequence $\set{\salg{A}_n}_{n\in \nset{N}}$ of $\sigma$-subalgebras of $\salg{A}$, it it a natural problem to establish conditions under which the conditional expectations converge almost everywhere.

In [] we analyzed the case for a  sequence to converge in $L^p$   and provided necessary an sufficient conditions for that to happen. We defined two  $\sigma$-subalgebras $\salg{A}_\mu$ and $\salg{A}^\perp$ and proved that we have convergence in $L^p$ if and only if $\salg{A}_\mu = \salg{A}^\perp$. We also established a chain of $\sigma$-subalgebras.

\begin{eqnarray*}
\underline{\salg{A}}\subset  \salg{A}_\mu \subset \salg{A}^\perp \subset \overline{\salg{A}}.
\end{eqnarray*} 

where $\underline{\mathfrak{A}}$ and $\overline{\mathfrak{A}}$ are the inferior and the superior limit of the  the $\sigma$-subalgebras  $\set{\mathfrak{A}_n}_{n\in \nset{N}}$. That is $
\underline{\mathfrak{A}}=\bigvee_{m=1}^\infty \bigcap_{n=m}^\infty \mathfrak{A}_n$ and $ \overline{\mathfrak{A}}=\bigcap_{m=1}^\infty \bigvee_{n=m}^\infty $. As a corollary we had convergence in $L^p$ in the special case  $\underline{\mathfrak{A}}= \overline{\mathfrak{A}}$ which was a result previously obtained by Fetter []. In that paper she pondered if, as in the case of monotone convergence of $\sigma$-algebras, we have also  convergence a.e.. However, in [] a negative answer was provided by establishing a counterexample.

In this paper we will analyze a  chain of families of  sets

\begin{eqnarray*}
\underline{\mathfrak{A}}\subset \mathfrak{A}_{w.a.e.} \subset \mathfrak{A}_\mu \subset \mathfrak{A}^\perp \subset \mathfrak{A}^{\perp a.e.}\subset \overline{\mathfrak{A}}.
\end{eqnarray*}

We will show that we have a.e. convergence if and only if $ \mathfrak{A}_{w.a.e.}= \mathfrak{A}^{\perp a.e.}$ provided that these two sets are $\sigma$-subalgebras.
\

\section{Previous-Works}

As usual, we will use the notation $\mathcal{E}\left(f|\mathfrak{A} \right)$ for the conditional expectation of $f$ given the $\sigma$-algebra $\mathfrak{A}$.     
Some well-known results for $a.e.$-convergence and $L^p$-convergence ($1\leq p<\infty$) are :

\begin{theorem}\label{teo:01}
(Martingales). If $\set{\mathfrak{A}_n}_{n\in \nset{N}}$ is monotone increasing sequence of $\sigma$-subalgebras of $\mathfrak{A}$, that is, $\mathfrak{A}_{n}\subset \mathfrak{A}_{n+1}$ for any $n\in \nset{N}$, then

\begin{eqnarray*}
\mathcal{E}\left(f|\mathfrak{A}_n \right)\underset{L^p}{ \overset{a.e.}{\longrightarrow}} \mathcal{E}\left(f\left|\bigvee_{n=1}^\infty \mathfrak{A}_n \right. \right),
\end{eqnarray*}

for every $f\in L^p(\mathfrak{A})$, where $\bigvee_{n=1}^\infty \mathfrak{A}_n$ stands for the minimum $\sigma$-algebra that contains $\bigcup_{n=1}^\infty \mathfrak{A}_n$. \big{(} Or if $\set{\mathfrak{A}_n}_{n\in \nset{N}}$ is monotone decreasing, that is, $\mathfrak{A}_n\supset \mathfrak{A}_{n+1}$, then $\mathcal{E}\left(f|\mathfrak{A}_n \right) \underset{L^p}{ \overset{a.e.}{\longrightarrow}} \mathcal{E}\left(f\left|\bigcap_{n=1}^\infty \mathfrak{A}_n \right.\right)$ \big{)}.

\end{theorem}

In [1], Boylan introduced a Hausdorff metric of the space of $\sigma$-algebras. It gives us a relationship between Cauchy sequences of $\sigma$-subalgebras and $L^p$-convergence  of conditional expectations.

\begin{theorem}\label{teo:03}
(Boylan, Equiconvergence). Let $\set{\mathfrak{A}_n}_{n\in \nset{N}}$ be a Cauchy sequence on the space of $\sigma$-algebras with the Hausdorff metric, that is,

\begin{eqnarray*}
d(\mathfrak{A}_n,\mathfrak{A}_m)=\sup_{A\in \mathfrak{A}_n}\left(\inf_{B\in \mathfrak{A}_m}\mu(A  \bigtriangleup B) \right)+ \sup_{B\in \mathfrak{A}_n}\left(\inf_{A\in \mathfrak{A}_m} \mu(A  \bigtriangleup B) \right),
\end{eqnarray*}

here $A\bigtriangleup B=(A\setminus B) \cup (B\setminus A)$. There is a subalgebra $\mathfrak{D}$ such that 

\begin{eqnarray*}
\lim_{n\to \infty}d(\mathfrak{A}_n,\mathfrak{D})=0,
\end{eqnarray*}

and

\begin{eqnarray*}
\mathcal{E}\left(f|\mathfrak{A}_n \right) \overset{L^p}{\longrightarrow} 
\mathcal{E}\left(f|\mathfrak{D} \right), 
\end{eqnarray*}

for every $f\in L^p(\mathfrak{A}), \ 1\leq p < \infty$ (see [2,3]).

\end{theorem}

Another approach was given by Fetter([]). She proved that if the lim sup of a  sequence of  $\sigma$-algebras coincides with the lim inf, then we have convergence in $L^p$. Indeed  

\begin{theorem}\label{teo:02}
(Fetter). If $\set{\mathfrak{A}_n}_{n\in \nset{N}}$ is such that $\underline{\mathfrak{A}}=\overline{\mathfrak{A}}$ where

\begin{eqnarray*}
\underline{\mathfrak{A}}=\bigvee_{m=1}^\infty \bigcap_{n=m}^\infty \mathfrak{A}_n & \mbox{and}& \overline{\mathfrak{A}}=\bigcap_{m=1}^\infty \bigvee_{n=m}^\infty \mathfrak{A}_n,
\end{eqnarray*}

then

\begin{eqnarray*}
\mathcal{E}\left(f|\mathfrak{A}_n \right)\underset{L^p}{\longrightarrow} \mathcal{E}\left(f\left|\overline{\mathfrak{A}} \right. \right),
\end{eqnarray*}

for every $f\in L^p(\mathfrak{A}), \ 1\leq p <\infty$ (see [1,4]).

\end{theorem}

Since the condition of the above theorem is fulfilled when we have monotone sequences of $\sigma$-algebras, Fetter´s result implies that of the martingale theorem in the case of $L^p$ convergence. However, we point out that in [4] it was proved that the condition $\underline{\mathfrak{A}}=\overline{\mathfrak{A}}$ does not imply convergence almost everywhere.

Given a sequence $\{ \mathfrak{A}_n \}_{n \in \nset{N}}$ of $\sigma$-algebras, two relevant $\sigma$-algebras were defined in []:

\begin{eqnarray*}
\mathfrak{A}_\mu=\set{A\in \mathfrak{A} \ : \ \exists A_n\in \mathfrak{A}_n, \hspace{0.1cm}\lim_{n\to \infty}\mu(A_n \bigtriangleup A)=0  },
\end{eqnarray*}

and if

\begin{eqnarray*}
W= \set{ g\in L^2( \mathfrak{A}) \ : \exists  A_{n_k} \in  \mathfrak{A}_{n_k}, \text{with }  \chi_{ A_{n_k}}\overset{weakly}{\longrightarrow} g} 
\end{eqnarray*}

$\mathfrak{A}_\perp$ was defined as the  minimum complete $\sigma$-algebra such that $g$  is $\mathfrak{A}_\perp$ measurable
for all $ g \in W$ . The importance of these $\sigma$-algebras is clear due to the following theorem:

\begin{theorem}\label{teo:04}
(Alonso-Brambila). Let $\set{\mathfrak{A}_n}_{n\in \nset{N}}$ be such that there is a $\sigma$-subalgebra $\mathfrak{A}_\infty$ with the property that $\mathfrak{A}_n \overset{\mu}{\to}\mathfrak{A}_\infty$ and $\mathfrak{A}_n \overset{\perp}{\to}\mathfrak{A}_\infty$. Then and only then

\begin{eqnarray*}
\mathcal{E}\left(f|\mathfrak{A}_n \right)\overset{L^p(\mathfrak{A})}{\longrightarrow} \mathcal{E}\left(f\left|\mathfrak{A}_\infty \right. \right),
\end{eqnarray*}

for every $f\in L^p(\mathfrak{A}), \ 1\leq p <\infty$.

\end{theorem}

Since these $\sigma$-algebras satisfy the relationship  

\begin{eqnarray*}
\underline{\mathfrak{A}}\subseteq \mathfrak{A}_\mu \subseteq \mathfrak{A}_\perp \subseteq \overline{\mathfrak{A}}.
\end{eqnarray*}

We have that  if $\underline{A}=\overline{A}$  Fetter´s theorem becomes  an immediate corollary.

Finally we point out that none of the theorems but theorem \ref{teo:01} deal with the problem of convergence a.e.

\section{A sigma-subalgebra}

In $[5]$ it was introduced the concept for a sequence of $\sigma$-subalgebras $\set{\mathfrak{A}_n}$ to $\mu$-approach  a $\sigma$-subalgebra $\mathfrak{D}$ if for each $D\in \mathfrak{D}$ there were $A_n\in \mathfrak{A}_n$ such that $\mu(A_n \bigtriangleup D)\to 0$. It was established that in such a case and only in that case we have for $f\in L^p(\mathfrak{D}) \hspace{0.1cm} (1\leq p<\infty)$ that

\begin{eqnarray*}
\mathcal{E}(f|\mathfrak{A}_n)\overset{L^p}{\longrightarrow} \mathcal{E}(f|\mathfrak{D})=f.
\end{eqnarray*}

It is clear that in order to have convergence $a.e.$ we need to introduce a more stronger concept for sets in $\set{\mathfrak{A}_n}$ to approach those of $\mathfrak{D}$.

As usual, in the following $A^c$ will stand for $\nset{X}\setminus A$, $\chi_A$ for the characteristic function of a set $A$, $A \bigtriangleup B$ for the symmetric difference of the sets $A$ and $B$, and $A=B \hspace{0.1cm} a.e$ for $\mu(A \bigtriangleup B)=0$.

We begin by first establishing the following lemma. 

\begin{lemma}\label{lem:08}
Let $\left( \nset{X},\mathfrak{A},\mu\right)$ a probability space with $\sigma$-algebra $\mathfrak{A}$ and measure $\mu$. If $\set{A_n}_{n\in \nset{N}}$ is a sequence of elements in $\mathfrak{A}$ and $A\in \mathfrak{A}$, the following statements are equivalent:

\begin{itemize}
\item[i)]$\chi_{A_n}\underset{a.e}{\longrightarrow}\chi_{A}$.\vspace{0.1cm}
\item[ii)]$\ds A=\bigcup_{N=1}^\infty \bigcap_{n\geq N}A_n=\bigcap_{N=1}^\infty \bigcup_{n\geq N}A_n \hspace{0.1cm}a.e.$\vspace{0.1cm}
\item[iii)]$\ds \lim_{N\to \infty}\mu\left(\bigcup_{n>N}(A_n \bigtriangleup A) \right)=0$.
\end{itemize}

\begin{proof}
Notice that $\chi_{A_n}\to \chi_A \hspace{0.1cm}a.e$ implies that, for all $x\in \nset{X}\setminus M$, where $M$ is a set of measure zero, there is an $N_x\in \nset{N}$ such that if $n>N_x$

\begin{eqnarray}\label{eq:01}
\abs{\chi_{A_n}(x)-\chi_A(x)}<\dfrac{1}{2}.
\end{eqnarray}

To prove $i)\to ii)$, we first take a look at the elements of $A\setminus M$. Since

\begin{eqnarray}
\abs{\chi_{A_n}(x)-\chi_A(x)}=\abs{\chi_{A_n^c}(x)-\chi_{A^c}(x)},
\end{eqnarray}

\eqref{eq:01} and  \eqref{eq:02} tell us that for $n>N_x$ and $x\in A\setminus M$; $\chi_{A_n^c}(x)<1/2$, and so $x\in A_n$. Therefore

\begin{eqnarray*}
A\setminus M \subset \bigcup_{N=1}^\infty \bigcap_{n>N}A_n.
\end{eqnarray*}

Using the same argument for $A^c$ we get

\begin{eqnarray*}
A^c\setminus M \subset \bigcup_{N=1}^\infty \bigcap_{n>N} A_n^c.
\end{eqnarray*}

Thus

\begin{eqnarray*}
A\cup M \supset \bigcap_{N=1}^\infty \bigcup_{n>N}A_n.
\end{eqnarray*}

Since

\begin{eqnarray*}
\bigcup_{N=1}^\infty \bigcap_{n>N}A_n \subset \bigcap_{N=1}^\infty \bigcup_{n>N}A_n,
\end{eqnarray*}

we have

\begin{eqnarray*}
\left(\bigcap_{N=1}^\infty \bigcup_{n>N}A_n \right)\cap M^c \subset A\cap M^c \subset \bigcup_{N=1}^\infty \bigcap_{n>N}A_n \subset \bigcap_{N=1}^\infty \bigcup_{n>N}A_n.
\end{eqnarray*}

So

\begin{eqnarray*}
\left( \bigcup_{N=1}^\infty \bigcap_{n>N}A_n \right)\cap M^c =\left( \bigcap_{N=1}^\infty \bigcup_{n>N}A_n \right)\cap M^c= A\cap M^c.
\end{eqnarray*}

Therefore $i\to ii)$.

We will prove now that $iii) \to i)$. Let

\begin{eqnarray*}
M=\bigcap_{N=1}^\infty \bigcup_{n>N}(A_n \bigtriangleup A).
\end{eqnarray*}

Then, by hypothesis $\mu(M)=0$. Now, as 

$$M^c=
\nset{X}\setminus M= \bigcup_{N=1}^\infty \bigcap_{n>N}(A_n \bigtriangleup A)^c.
$$

if $x\in \nset{X}\setminus M$ there is an $N\in \nset{N}$ such that

\begin{eqnarray*}
x\in \bigcap_{n>N} (A_n \bigtriangleup A )^c,
\end{eqnarray*} 

and hence $x\in (A_n \bigtriangleup A)^c$ for all $n>N$. That is

\begin{eqnarray*}
0=\chi_{A_n \bigtriangleup A}(x)=\abs{\chi_A(x)-\chi_{A_n}(x)}<\epsilon, & \mbox{for }n>N.
\end{eqnarray*}

Finally, to prove that $ii)$ implies $iii)$, we notice that:

\begin{eqnarray*}
\mu\left( \bigcup_{n>N}(A_n \bigtriangleup A)\right)=\mu\left( \left(\bigcup_{n>N}A_n \right)\cap A^c \right)+\mu\left(\left(\bigcup_{n>N}A_n^c \right)\cap A \right).
\end{eqnarray*}

Since by hypothesis

\begin{align*}
0=\mu\left(A\bigtriangleup \bigcup_{N=1}^\infty \bigcap_{n>N}A_n \right)\geq \mu\left(A\setminus \bigcup_{N=1}^\infty \bigcap_{n>N}A_n \right)=\mu\left( A\cap \bigcap_{N=1}^\infty\bigcup_{n>N}A_n^c\right)=\lim_{N\to \infty}\mu\left( A\cap \bigcup_{n>N}A_n^c\right),
\end{align*}

and

\begin{eqnarray*}
0=\mu\left(A\bigtriangleup \bigcap_{N=1}^\infty \bigcup_{n>N}A_n \right)\geq \mu\left(\bigcap_{N=1}^\infty \bigcup_{n>N}A_n\setminus A \right)=\lim_{N\to\infty}\mu\left( \bigcup_{n>N}A_n\cap A^c \right),
\end{eqnarray*}

we get

\begin{eqnarray*}
\lim_{N\to \infty}\mu\left(\bigcup_{n>N}(A_n\bigtriangleup A) \right)=0
\end{eqnarray*}

\end{proof}
\end{lemma}

Given a sequence $\set{\mathfrak{A}_n}_{n\in \nset{N}}$ of $\sigma$-subalgebras of $\mathfrak{A}$ we will define the family of sets  $\mathcal{F}$ as

\begin{eqnarray*}
\mathcal{F}=\set{A\in \mathfrak{A} \ : \ A=\bigcup_{N=1}^\infty \bigcap_{n>N}A_n=\bigcap_{N=1}^\infty \bigcup_{n>N}A_n \hspace{0.1cm} a.e. \mbox{ for a sequence }\set{A_n} \mbox{ such that } A_n\in \mathfrak{A}_n} .
\end{eqnarray*}

We will prove that $\mathcal{F}$ is a $\sigma$-algebra. To do this we will need the following lemma.

\begin{lemma}\label{lem:09}
If for any $r>0$ there is a sequence $\set{A_n^r \ : \ A_n^r\in \mathfrak{A}_n}$ and $N_r\in \nset{N}$ such that

\begin{eqnarray*}
\mu\left(\bigcup_{n>N_r}(A_n^r \bigtriangleup B) \right)<r,
\end{eqnarray*}

then $B\in \mathcal{F}$.

\begin{proof}
Let $S_m=1/2^m$ for $m\in \nset{N}$ and $\set{Q_m}_{m\in \nset{N}}$ a strictly increasing sequence such that $Q_m\geq N_{S_m}$. As $Q_n\geq N_{S_m}$ we have

\begin{eqnarray*}
\mu\left(\bigcup_{n>Q_m} (A_n^{S_m}\bigtriangleup B ) \right)\leq \mu\left(\bigcup_{n>N_{S_m}}( A_n^{S_m} \bigtriangleup B )\right)<S_m.
\end{eqnarray*}

For $n>Q_1$ define the sequence $\set{A_n}$ by:

\begin{eqnarray*}
A_n=A_n^{S_k}&\mbox{if}& Q_k<n\leq Q_{k+1},
\end{eqnarray*}

since $A_n^{S_k}\in \mathfrak{A}_n$ so is $A_n$. We have

\begin{align*}
\mu\left( \bigcup_{n>Q_m}(A_n \bigtriangleup B) \right)=& \mu\left(\bigcup_{k=m}^\infty \,\bigcup_{Q_k<n\leq Q_{k+1}}(A_n^{S_k} \bigtriangleup B) \right)\leq \sum_{k=m}^\infty \mu\left(\bigcup_{Q_k<n\leq Q_{k+1}}(A_n^{S_k}\bigtriangleup B )\right)\\
\leq & \sum_{k=m}^\infty \mu\left( \bigcup_{n>Q_k}(A_n^{S_k}\bigtriangleup B )\right)<\sum_{k=m}^\infty S_k=\dfrac{1}{2^{m-1}},
\end{align*}

and therefore

\begin{eqnarray*}
\lim_{N\to \infty}\mu\left( \bigcup_{n>N}(A_n\bigtriangleup B) \right)=0.
\end{eqnarray*}

The proof follows from lemma 3.1

\end{proof}

We can now prove that $\mathcal{F}$ is a $\sigma$-algebra.

\begin{proposition}\label{prop:01}
$\mathcal{F}$ is a $\sigma$-subalgebra of $\mathcal{A}$

\begin{proof}
It is clear that $\emptyset,\ \nset{X}$ are elements of $\mathcal{F}$ since we can take $A_n=\emptyset$ for all $n\in \nset{N}$ or $A_n=\nset{X}$ for all $n\in \nset{N}$.

If $A\in \mathcal{F}$, by Proposition 3.1 there are $\set{A_n\in \mathfrak{A}_n}$ such that

\begin{eqnarray*}
\chi_{A_n}\longrightarrow \chi_{A} \hspace{0.1cm}a.e.
\end{eqnarray*}

Since $A_n^c \in \mathfrak{A}_n$ and

\begin{eqnarray*}
\chi_{A_n^c}=1-\chi_{A_n}\overset{a.e.}{\longrightarrow} 1-\chi_A=\chi_{A^c},
\end{eqnarray*}

we have that $A^c\in \mathcal{F}$.

Let $A,B\in \mathcal{F}$. By proposition 3.1 there are two sequences $\set{A_n}, \ \set{B_n}$, with  $A_n,B_n\in \mathfrak{A}_n$ for all $n\in \nset{N}$,  such that $\chi_{A_n}\to \chi_{A} \hspace{0.1cm}a.e.$ and $\chi_{B_n}\to \chi_B \hspace{0.1cm}a.e.$ Since $A_n\cap B_n\in \mathfrak{A}_n$ and $\chi_{A_n\cap B_n}=\chi_{A_n}\chi_{B_n}\underset{a.e.}{\longrightarrow} \chi_{A}\chi_{B}=\chi_{A\cap B}$ we have that $A\cap B \in \mathcal{F}$. Thus $\mathcal{F}$ is an algebra.

To show that is is a $\sigma$-algebra, let $B=\bigcup_{k=1}^\infty D_k$ with $D_k\in \mathcal{F}$. Since $\mathfrak{F}$ is an algebra, we can consider the sets  $\{D_k\}$ disjoint.

Define $E_M=\bigcup_{k=1}^M D_k$ and let $r>0$. Since the sets $D_k$ are disjoint, there is an $M_1\in \nset{N}$ such that

\begin{eqnarray*}
\mu\left( \bigcup_{k>M_1}^\infty D_k \right)<\dfrac{r}{2}.
\end{eqnarray*}

On the other hand, since ${E_M}_1\in \mathcal{F}$. There is a sequence $\set{A_n\in \mathfrak{A}_n}$ such that

\begin{eqnarray*}
\lim_{N\to \infty}\mu\left( \bigcup_{n>N} (A_n \bigtriangleup E_{M_1}) \right)=0.
\end{eqnarray*}

Therefore, there is an $N_1\in \nset{N}$ such that $\mu\left( \bigcup_{n>N_1}(A_n\bigtriangleup E_{M_1}) \right)<r/2$.

As $E_{M_1}\subset B$, $A_n\setminus B \subset A_n\setminus E_{M_1}$. So

\begin{eqnarray*}
\mu\left(\bigcup_{n>N_1}A_n\setminus B \right)\leq \mu\left(\bigcup_{n>N_1}(A_n\setminus E_{M_1}) \right)\leq \mu\left(\bigcup_{n>N_1}(A_n\bigtriangleup E_{M_1}) \right)<\dfrac{r}{2}.
\end{eqnarray*}

Also, as $B$ is the union of the disjoint sets $E_{M_1}$ and $\bigcup_{k>M_1}^\infty D_k$ we have

\begin{align*}
\mu\left(\bigcup_{n>N_1}(B\setminus A_n) \right)=& \mu\left( \bigcup_{n>N_1}(E_{M_1}\setminus A_n) \right)+\mu\left( \bigcup_{n>N_1}\left(\bigcup_{k>M_1}^\infty D_k \right)\setminus A_n \right)\\
\leq & \mu\left( \bigcup_{k>N_1}(E_{M_1} \bigtriangleup A_n) \right)+\mu\left( \bigcup_{k>M_1}^\infty D_k \right)<\dfrac{r}{2}+\dfrac{r}{2}.
\end{align*}

Therefore $\mu\left( \bigcup_{n>N_1}(A_n \bigtriangleup B) \right)<r$.

So, by lemma 3.2 $B$ is in $\mathcal{F}$.
\end{proof}

\end{proposition}

\end{lemma}

Given a sequence of $\sigma$-subalgebras of $\mathfrak{A}$ in $[5]$ we found for the $\sigma$-algebra $\mathfrak{A}_\mu$ defined by

\begin{eqnarray*}
\mathfrak{A}_\mu=\set{A\in \mathfrak{A} \ : \ \exists A_n\in \mathfrak{A}_n, \hspace{0.1cm}\lim_{n\to \infty}\mu(A_n \bigtriangleup A)=0  },
\end{eqnarray*}

the relationships

\begin{eqnarray*}
\bigvee_{m=1}^\infty \bigcap_{n=m}^\infty \mathfrak{A}_n \equiv \underline{\mathfrak{A}} \hspace{0.1cm}\subset \hspace{0.1cm} \mathfrak{A}_\mu  \hspace{0.1cm} \subset \hspace{0.1cm} \overline{\mathfrak{A}}\equiv \bigcap_{m=1}^\infty \bigvee_{n=m}^\infty \mathfrak{A}_n.
\end{eqnarray*}

It is clear by its definition that $\mathcal{F}\subset \mathfrak{A}_\mu$. Now, let $A\in \bigcap_{n=N}^\infty \mathfrak{A}_n$, take $A_n=A$ for $n>N$, then trivially $\chi_{A_n}\underset{a.e.}{\longrightarrow}\chi_A$ and thus $A\in \mathcal{F}$. Since $\mathcal{F}$ is a $\sigma$-algebra, we have

\begin{eqnarray*}
\mathcal{F}\supset \bigvee_{N=1}^\infty \bigcap_{n=N}^\infty \mathfrak{A}_n.
\end{eqnarray*}

Therefore

\begin{eqnarray*}
\underline{\mathfrak{A}}\subset \mathcal{F}\subset \mathfrak{A}_\mu \subset \overline{\mathfrak{A}}.
\end{eqnarray*}

What we have proven is that given a sequence of $\sigma$-algebras $\set{\mathfrak{A}_n}$, if $A$ is a set $\mathcal{F}$ measurable, there are $f_n$ functions $\mathfrak{A}_n$ measurable such that $f_n\to \chi_A \hspace{0.1cm}a.e.$ However that does not imply that $\mathcal{E}(\chi_A|\mathfrak{A}_n)\to \chi_A \hspace{0.1cm}a.e.$ even when the conditional expectations do so in $L^p \hspace{0.1cm}(1\leq p< \infty)$.

In the example shown in $[4]$ we have the property that $\underline{\mathfrak{A}}=\overline{\mathfrak{A}}$ and hence $\underline{\mathfrak{A}}=\mathcal{F}$, and a Borel set for which the conditional expectations do not converge $a.e.$ to its characteristic function. Let us show an easier example for which the conditional expectations for a characteristic function in $\mathcal{F}$ do not converge $a.e.$ but they do so in $L^p \hspace{0.1cm}(1\leq p<\infty)$.

Let $\nset{X}=[0,1]$ and $\mathfrak{A}$ the Lebesgue measurable sets. We are going to define a sequence $\set{A_n}$ that lie in $\left[ \frac{1}{2},1\right)$ and such that $\chi_{A_n}\to \chi_{\left[\frac{1}{2},1\right)} \hspace{0.1cm}a.e.$ $\set{J_k}$ will be a sequence that lies in $\left[\frac{1}{2},1\right)\setminus A_n$ but with measure half its size

Specifically

\begin{eqnarray*}
A_n=\left[\dfrac{1}{2},1-\dfrac{1}{2^n} \right), & J_n=\left[1-\dfrac{1}{2^{n+1}},1 \right), & n\in \nset{N}.
\end{eqnarray*}

On $\left[0,\frac{1}{2}\right)$ define a sequence of intervals that go back and forth in the interval and in each turn diminishing its size

\begin{eqnarray*}
I_{n,k}=\left[\dfrac{k}{4\cdot 2^n}, \dfrac{k+1}{4\cdot2^n} \right)&\mbox{ with }0\leq k<2\cdot 2^n, & n\in \nset{N}.
\end{eqnarray*}

We have that for each $n\in \nset{N}$

\begin{eqnarray*}
\bigcup_{k=0}^{2\cdot 2^n-1}I_{n,k}=\left[ 0,\dfrac{1}{2}\right).
\end{eqnarray*}

Let $B_{n,k}=I_{n,k}\cup J_n$ and $C_{n,k}=(A_n\cup B_{n,k})^c=A_n^c \cap B_{n,k}^c$.

Consider the sequence of $\sigma$-subalgebras

\begin{eqnarray*}
\mathfrak{A}_{n,k}=\set{\emptyset, \nset{X},A_n, B_{n,k},C_{n,k}},
\end{eqnarray*}

Since

\begin{eqnarray*}
\mathcal{E}\left(\chi_{\left[\frac{1}{2},1\right)}\big{|} \mathfrak{A}_{n,k} \right)=\dfrac{
\langle\chi_{\left[\frac{1}{2},1\right)},\chi_{A_n} \rangle}{\mu(A_n)}\chi_{A_n}+\dfrac{\langle\chi_{\left[\frac{1}{2},1\right)},\chi_{B_{n,k}} \rangle}{\mu(B_{n,k})}\chi_{B_{n,k}}+  \dfrac{\langle\chi_{\left[\frac{1}{2},1\right)},\chi_{C_{n,k}} \rangle}{\mu(C_{n,k})}\chi_{C_{n,k}},
\end{eqnarray*}

we have

\begin{align*}
\mathcal{E}\left(\chi_{\left[\frac{1}{2},1\right)}\big{|}\mathfrak{A}_{n,k} \right)=&\chi_{A_n}+\dfrac{\dfrac{1}{2^{n+1}}}{\dfrac{3}{2}\dfrac{1}{2^{n+1}}}\chi_{B_{n,k}}+ \dfrac{\dfrac{1}{2^{n+1}}}{\dfrac{1}{2}+\dfrac{1}{2^{n+2}}}\chi_{C_{n,k}}\\
=&\chi_{A_n}+\dfrac{2}{3}\chi_{B_{n,k}}+\dfrac{2}{1+2^{n+1}}\chi_{C_{n,k}}\\
=&\chi_{A_n}+\dfrac{2}{3}\chi_{I_{n,k}}+\dfrac{2}{3}\chi_{J_n}+\dfrac{2}{1+2^{n+1}}\chi_{C_{n,k}}.
\end{align*}

Since $\chi_{A_n}\to \chi_{\left[\frac{1}{2},1\right)} \hspace{0.1cm}a.e.$, $\dfrac{2}{3}\chi_{J_n}\to 0 \hspace{0.1cm}a.e.$ and $\dfrac{2}{1+2^{n+1}}\chi_{C_{n,k}}\to 0 \hspace{0.1cm}a.e.$ but $\dfrac{2}{3}\chi_{I_{n,k}}$ does not convergence $a.e.$, we have that $\left[\frac{1}{2},1\right)\in \mathcal{F}$ but $\mathcal{E}\left(\chi_{\left[\frac{1}{2},1\right)}\big{|}\mathfrak{A}_{n,k}\right)\nrightarrow \chi_{\left[\frac{1}{2},1\right)} \hspace{0.1cm}a.e.$ 

\section{Necessary and sufficient conditions for almost everywhere convergence.}

\subsection{On almost everywhere convergence of characteristics functions}

Let $\mathfrak{B}$ be a $\sigma$-subalgebra of $\mathfrak{A}$.

\begin{definition}\label{def:01}
Define the seminorm $\norm{ \ \cdot \ }_{\mathfrak{B}}$ for $f\in L^\infty(d\mu)$ as

\begin{eqnarray*}
\norm{f}_{\mathfrak{B}}=\norm{\mathcal{E}\left(f|\mathfrak{B} \right)}_\infty.
\end{eqnarray*}

\end{definition}

In the following we are going to use the notation   $\norm{A}_{\mathfrak{B}}=\norm{\chi_A}_{\mathfrak{B}}$ for any set $A \in \mathfrak{A}$.  We have the property

\begin{lemma}\label{lem:01}
Let $A$ be a measurable set in $\mathfrak{A}$ and $\chi_A$ its characteristic function. Then

\begin{eqnarray*}
\norm{A}_{\mathfrak{B}}=\norm{\chi_A}_{\mathfrak{B}}=\sup_{\tiny \begin{matrix}
B\in \mathfrak{B}\\
\mu(B)>0
\end{matrix}} \dfrac{\mu(A\cap B)}{\mu(B)}.
\end{eqnarray*}

\begin{proof}
Let $B\in \mathfrak{B}$. We have then

\begin{align*}
\mu(A\cap B)&=\int \chi_{A\cap B}d\mu=\int \chi_A \chi_B d\mu \\
&= \int \mathcal{E}(\chi_A|\mathfrak{B})\chi_Bd\mu \leq \norm{\mathcal{E}(\chi_A|\mathfrak{B})}_\infty \mu(B).
\end{align*}

Therefore, for any $B\in \mathfrak{B}$ with $\mu(B)>0$

\begin{eqnarray*}
\dfrac{\mu(A\cap B)}{\mu(B)}\leq  \norm{\mathcal{E}(\chi_A|\mathfrak{B})}_\infty.
\end{eqnarray*}

Since the case $\norm{\mathcal{E}(\chi_A|\mathfrak{B})}_\infty=0$ is trivial, take any $\varepsilon>0$ such that

\begin{eqnarray*}
\norm{\mathcal{E}(\chi_A|\mathfrak{B})}_\infty-\varepsilon>0.
\end{eqnarray*}

By definition the set

\begin{eqnarray*}
C=\set{x  \ : \ \mathcal{E}(\chi_A|\mathfrak{B})>   \norm{\mathcal{E}(\chi_A|\mathfrak{B})}_\infty-\varepsilon }
\end{eqnarray*}

is $\mathfrak{B}$-measurable and

\begin{eqnarray*}
\int_C \mathcal{E}(\chi_A|\mathfrak{B})d\mu \geq \left( \norm{\mathcal{E}(\chi_A|\mathfrak{B})}_\infty-\varepsilon\right)\mu(C).
\end{eqnarray*}

as

\begin{eqnarray*}
\int_C \mathcal{E}(\chi_A|\mathfrak{B})d\mu= \int \chi_C \mathcal{E}(\chi_A | \mathfrak{B})d\mu= \int \chi_C \chi_A d\mu= \mu(C\cap A),
\end{eqnarray*}

we have then

\begin{eqnarray*}
\dfrac{\mu(C\cap A)}{\mu(C)}\geq \norm{\mathcal{E}(\chi_A|\mathfrak{B})}_\infty-\varepsilon.
\end{eqnarray*}

\end{proof}

\end{lemma}

\begin{definition}\label{def:02}
We say that $A\in \mathfrak{A}$ is uniformly covered by the sequence  of $\sigma$-subalgebras $\set{\mathfrak{A}_n}_{n\in \nset{N}}$ if there is a sequence $\set{A_n\in \mathfrak{A}_n}_{n\in \nset{N}}$ such that

\begin{itemize}
\item[i)]$\ds A=\bigcup_{N=1}^\infty \bigcap_{n\geq N}A_n= \bigcap_{N=1}^\infty \bigcup_{n\geq N}A_n$. \vspace{0.3cm}
\item[ii)] $\ds \norm{\chi_{A\setminus A_n}}_{\mathfrak{A}_n}\to 0$   \hspace{0.1cm} as \hspace{0.1cm} $n\to \infty$.
\end{itemize}

\end{definition}

Next lemma shows that actually, we can relax a little bit condition $ii)$.

\begin{lemma}\label{lem:04}

$A$ is uniformly covered by $\set{\mathfrak{A}_n}_{n\in \nset{N}}$ if and only if for any $r>0$ there is a sequence $\set{A_n^r}_{n\in \nset{N}}$, $A_n^r\in \mathfrak{A}_n$ and $M\in \nset{N}$ such that

\begin{itemize}
\item[i)]$\ds A=\bigcup_{N=1}^\infty \bigcap_{n>N}A_n^r=\bigcap_{N=1}^\infty \bigcup_{n>N}A_n^r$ \vspace{0.4cm}.

\item[ii)] $\norm{A\setminus A^r}_{\mathfrak{A}_n}<r$ \hspace{0.2cm} if \hspace{0.2cm} $n>M_r$.
\end{itemize}

\begin{proof}
By lemma 3.1 the condition $i)$ means that 

\begin{eqnarray*}
\lim_{N\to \infty}\mu \left( \bigcup_{n>N}(A  \bigtriangleup A_n^r)\right)=0.
\end{eqnarray*}

Thus, for $r=1$ let $N_1>M_1$ and such that

\begin{eqnarray*}
\mu \left(\bigcup_{n>N_1} (A\bigtriangleup A_n^r\right)<1
\end{eqnarray*}.

In general let $r_k=2^{-k}$ and $\tilde{N}_k$ such that

\begin{eqnarray*}
\mu\left( \bigcup_{n>\tilde{N}_k}(A \bigtriangleup A_n^k) \right) < r_k,
\end{eqnarray*}

and ${N}^\prime_k$ such that 

\begin{eqnarray*}
\norm{A\setminus A_n^{r_k}}_{\mathfrak{A}_n}<r_k & \mbox{ for }n>N^\prime_k.
\end{eqnarray*}

It is clear that we can define a strictly increasing sequence $N_k$ such that $N_k > \max{(\tilde{N}_k,N^\prime_k)}$ 

Define $A_n \in \mathfrak{A}_n$ as $A_n^{r_k}$  if $N_k\leq n \leq N_{k+1}$. Then

\begin{eqnarray*}
\norm{A\setminus A_n}_{\mathfrak{A}_n}=\norm{A\setminus A_n^{r_k}}_{\mathfrak{A}_n}<\dfrac{1}{2^k},
\end{eqnarray*}
we also have

\begin{align*}
\mu\left(\bigcup_{n\geq \widetilde{N}_{k'}} (A\bigtriangleup A_n)\right)=\mu\left(\bigcup_{k=k'}^\infty \bigcup_{\widetilde{N}_{k}\leq n <\widetilde{N}_{k+1}}(A\bigtriangleup A_n) \right)\leq \sum_{k=k'}^\infty \mu\left( \bigcup_{\widetilde{N}_k\leq n<\widetilde{N}_{k+1}  } A \bigtriangleup A_n^{r_k} \right)<\sum_{k=k'}^\infty \dfrac{1}{2^k}=\dfrac{1}{2^{k'-1}}.
\end{align*}

Therefore since $\mu\left(\bigcup_{n>N}A\bigtriangleup A_n \right)$ is monotone

\begin{eqnarray}
\lim_{N\to \infty} \mu\left(\bigcup_{n\geq N}A\bigtriangleup A_n \right)=0.
\end{eqnarray}

and thus finally

\begin{eqnarray*}
\lim_{n\to \infty}\norm{A\setminus A_n}_{\mathfrak{A}_n}=0.
\end{eqnarray*}

\end{proof}

\end{lemma}

We are going to say that a sequence of sets $\{A_n \in \mathfrak{A}_n\}$ uniformly covers a set $A \in \mathfrak{A}$ if conditions $i)$ and $ii)$ of definition XXX are satisfied. A couple of results regarding uniform covering are:

\begin{lemma}\label{lem:02}
 If $\{A_n \in \mathfrak{A}_n\}$ uniformly covers a set $A \in \mathfrak{A}$  and
  $\{{A^\prime}_n \in \mathfrak{A}_n\}$ is such that $
  A_n \subset  {A^\prime}_n$ for all $n \in \mathcal{N}$ and $\chi_{{A^\prime}_n} \longrightarrow  \chi_{A} \hspace{0.1cm}a.e. $
  then it uniformly covers $A$.
  \begin{proof}
  The proof is trivial since for $0 \leq f \leq g$, and $\mathfrak{B}$ $\sigma$-subalgebra
  $\mathcal{E}(f|\mathfrak{B})\leq \mathcal{E}(g|\mathfrak{B}) $ and so, as $A\setminus {A^\prime}_n \subset  A\setminus A_n$,
  $$\mathcal{E}(\chi_{A\setminus A^\prime_n}|\mathfrak{A}_n) \leq 
  \mathcal{E}(\chi_{A\setminus A_n}|\mathfrak{A}_n)
 $$
  \end{proof}
\end{lemma}
\begin{lemma}\label{lem:02}
If $A$ and $B$ are sets in $\mathfrak{A}$ uniformly covered by $\set{\mathfrak{A}_n}_{n\in \nset{N}}$. Then so are $A\bigcup B$ and $A \bigcap B$.

\begin{proof}
Let $A_n \in \mathfrak{A}_n$ and $B_n  \in \mathfrak{A}_n$ sequences of sets that uniformly cover $A$ and $B$ respectively. Consider first the sequence $C_n=A_n \bigcap B_n \in \mathfrak{A}_n$. Since by property $i)$ of the definition of uniform covering we have that $\chi_{A_n}\overset{a.e.}{\longrightarrow}\chi_A $ and $\chi_{B_n}\overset{a.e.}{\longrightarrow}\chi_B$:
$$
\chi_{A_n \bigcap B_n}=\chi_{A_n}\chi_{B_n}\overset{a.e.}{\longrightarrow}\chi_{A}\chi_{B}=\chi_{A \bigcap B}
$$

We have also that

$$
A \cap B \setminus (A_n \cap B_n)=A \cap B \cap   ({A^c}_n \cup {B^c}_n) \subset  (A \cap {A^c}_n) \cup (B \cap {B^c}_n)
$$
Therefore
$$
  \norm{A \cap B \setminus (A_n \cap B_n)}_{\mathfrak{A}_n} \leq  \norm{(A \cap {A^c}_n) \cup (B \cap {B^c}_n)}_{\mathfrak{A}_n}\leq  \norm{A \cap {A^c}_n}_{\mathfrak{A}_n}+ \norm{B \cap {B^c}_n}_{\mathfrak{A}_n} \longrightarrow 0 a.e.
$$
 For the case of the union of two sets  we have
 $$
\chi_{A_n \bigcup B_n}=\chi_{A_n} + \chi_{B_n}-\chi_{A_n}\chi_{B_n}\overset{a.e.}{\longrightarrow}\chi_{A}+\chi_{B} -\chi_{A}\chi_{B}=\chi_{A \bigcup B}
$$

and
$$
(A \cup B) \setminus (A_n \cup B_n)=(A \cup B) \cap   ({A^c}_n \cap {B^c}_n) \subset  (A \cap {A^c}_n) \cup (B \cap {B^c}_n)
$$

\end{proof}
\end{lemma}

\begin{lemma}\label{lem:02}
If $A$ is uniformly covered by $\set{\mathfrak{A}_n}_{n\in \nset{N}}$. Then

\begin{eqnarray*}
\varlimsup_{n\to \infty}\left(\mathcal{E}(\chi_A|\mathfrak{A}_n )\right)(x)\leq \chi_A(x) \ a.e.
\end{eqnarray*}

\begin{proof}
Let $\set{A_n}$ be a sequence with the properties  of \textbf{Definition \ref{def:02}}. Then

\begin{align*}
\mathcal{E}(\chi_A|\mathfrak{A}_n)&=\mathcal{E}(\chi_A \chi_{A_n}+\chi_A \chi_{A_n^c}|\mathfrak{A}_n)\\
&= \chi_{A_n}\mathcal{E}(\chi_{A}|\mathfrak{A}_n)+\mathcal{E}(\chi_{A\setminus A_n}|\mathfrak{A}_n)\\
&\leq \chi_{A_n}+\norm{\chi_{A\setminus A_n}}_{\mathfrak{A}_n}.
\end{align*}

as $A$ is uniformly covered by the sequence $\set{A_n}$ then $i)$ of \textbf{Definition \ref{def:02}} implies that

\begin{eqnarray*}
\chi_{A_n}\overset{a.e.}{\longrightarrow}\chi_A.
\end{eqnarray*}

So

\begin{eqnarray*}
\varlimsup_{n\to \infty}\mathcal{E}(\chi_A|\mathfrak{A}_n)\leq \chi_A \ a.e.
\end{eqnarray*}

\end{proof}

\end{lemma}

xxx222

In view of the proof we can actually relax somewhat the condition of uniformly covering to get a similar result.

\begin{lemma}\label{def:02}
Let $\set{\mathfrak{A}_n}_{n\in \nset{N}}$ be a sequence of $\sigma$-subalgebras . If  $A\in \mathfrak{A}$ is such that there is a sequence $\set{A_n\in \mathfrak{A}_n}_{n\in \nset{N}}$ satisfying

\begin{itemize}
\item[i)]$\ds A=\bigcup_{N=1}^\infty \bigcap_{n\geq N}A_n= \bigcap_{N=1}^\infty \bigcup_{n\geq N}A_n$. \vspace{0.3cm}
\item[ii)] $\mathcal{E}(\chi_{A\setminus A_n}|\mathfrak{A}_n) \to 0$ a.e   \hspace{0.1cm} as \hspace{0.1cm} $n\to \infty$.
\end{itemize}
Then 
\begin{eqnarray*}
\varlimsup_{n\to \infty}\left(\mathcal{E}(\chi_A|\mathfrak{A}_n )\right)(x)\leq \chi_A(x) \ a.e.
\end{eqnarray*}
\begin{proof}
The proof is exactly the same as the one of the above lemma.

\begin{align*}
\mathcal{E}(\chi_A|\mathfrak{A}_n)&=\mathcal{E}(\chi_A \chi_{A_n}+\chi_A \chi_{A_n^c}|\mathfrak{A}_n)\\
&= \chi_{A_n}\mathcal{E}(\chi_{A}|\mathfrak{A}_n)+\mathcal{E}(\chi_{A\setminus A_n}|\mathfrak{A}_n)\\
&\leq \chi_{A_n}+\mathcal{E}(\chi_{A\setminus A_n}|\mathfrak{A}_n).
\end{align*}

\end{proof}

\end{lemma}

Notice that if A is uniformly covered by $\{\mathfrak{A}_n\}$ the conditions of lemma 4.8xx are satisfied.

The interesting case occurs when both $A$ and $A^c$ are uniformly covered by $\set{\mathfrak{A}_n}$ or satisfy the conditions of the above lemma.

\begin{lemma}\label{lem:03}
If $A\in \mathfrak{A}$ is such that $A$ and $A^c$ are uniformly covered by $\set{\mathfrak{A}_n}$ (or satisy the conditions of lemma 4.9x) then

\begin{eqnarray*}
\mathcal{E}(\chi_A|\mathfrak{A}_n)\overset{a.e.}{\longrightarrow}\chi_A.
\end{eqnarray*}

\begin{proof}
\begin{align*}
\overline{\lim_{n\to \infty}}\mathcal{E}(\chi_{A^c}|\mathfrak{A}_n)&= \varlimsup_{n\to \infty}\mathcal{E}(1-\chi_A|\mathfrak{A}_n)\\
&=\varlimsup_{n\to\infty} \left(1-\mathcal{E}(\chi_A|\mathfrak{A}_n) \right)\\
&=1-\varliminf_{n\to\infty} \mathcal{E}(\chi_A|\mathfrak{A}_n).
\end{align*}

Since $A^c$ is uniformly covered

\begin{eqnarray*}
\chi_{A^c}=1-\chi_{A}\geq 1-\varliminf_{n\to \infty}\mathcal{E}(\chi_A|\mathfrak{A}_n),
\end{eqnarray*}

and so

\begin{eqnarray*}
\varliminf_{n\to \infty}\mathcal{E}(\chi_A|\mathfrak{A}_n)\geq \chi_A \ a.e.
\end{eqnarray*}

finally, as $A$ is uniformly covered

\begin{eqnarray*}
\varlimsup_{n\to\infty}\mathcal{E}(\chi_A|\mathfrak{A}_n)\leq \chi_A \leq \varliminf_{n\to\infty} \mathcal{E}(\chi_A|\mathfrak{A}_n) \ a.e.
\end{eqnarray*}

And then

\begin{eqnarray*}
\mathcal{E}(\chi_A|\mathfrak{A}_n)\overset{a.e.}{\longrightarrow}\chi_A.
\end{eqnarray*}

\end{proof}

\end{lemma}

Now the necessary lemma for $a.e.$ convergence of characteristic functions.

\begin{lemma}\label{lem:05}
Let $A\in \mathfrak{A}$ and $\set{\mathfrak{A}_n}_{n\in \nset{N}}$ $\sigma$-subalgebras such that

\begin{eqnarray*}
\mathcal{E}(\chi_A|\mathfrak{A}_n)\overset{a.e.}{\longrightarrow} \chi_A.
\end{eqnarray*}

Then $A$ and $A^c$ are uniformly covered by $\set{\mathfrak{A}_n}_{n\in \nset{N}}$.

\begin{proof}
Let $0<r<1$. Define $A_n\in \mathfrak{A}_n$ as

\begin{eqnarray*}
A_n=\set{x\in \nset{X} \ : \ \mathcal{E}(\chi_A|\mathfrak{A}_n)(x)\geq r},
\end{eqnarray*}

since $\mathcal{E}(\chi_A|\mathfrak{A}_n)\overset{a.e.}{\to} \chi_A$. We have that

\begin{eqnarray*}
\chi_{A_n}\overset{a.e.}{\longrightarrow}\chi_A,
\end{eqnarray*}

and therefore

\begin{eqnarray*}
A=\bigcap_{N=1}^\infty \bigcup_{n>N}A_n=\bigcup_{N=1}^\infty \bigcap_{n>N}A_n.
\end{eqnarray*}

It is also clear that

\begin{align*}
0\leq & \mathcal{E}(\chi_{A\setminus A_n}|\mathfrak{A}_n)=\mathcal{E}(\chi_A \chi_{A_n^c}|\mathfrak{A}_n)\\
=& \chi_{A_n^c}\mathcal{E}(\chi_A|\mathfrak{A}_n)< \chi_{A_n^c}r\leq r.
\end{align*}

Thus $\norm{\mathcal{E}(\chi_{A\setminus A_n}|\mathfrak{A}_n)}_\infty \leq r$ and by the above lemma $A$ is uniformly covered. 

Since $
\mathcal{E}(\chi_A|\mathfrak{A}_n)\overset{a.e.}{\longrightarrow}\chi_A$ implies that $
\mathcal{E}(\chi_{A^c}|\mathfrak{A}_n)\overset{a.e.}{\longrightarrow}\chi_{A^c}
$, we have that $A^c$ is also uniformly covered.

\end{proof}

\end{lemma}

From \textbf{Lemma \ref{lem:03}} and \textbf{Lemma \ref{lem:05}} we have the following theorem.

\begin{theorem}\label{teo:05}
Let $\set{\mathfrak{A}_n}_{n\in \nset{N}}$ be a sequence of $\sigma$-algebras and let $A\in \mathfrak{A}$. Then

\begin{eqnarray*}
\mathcal{E}(\chi_A|\mathfrak{A}_n)\overset{a.e.}{\longrightarrow}\chi_A,
\end{eqnarray*}

if and only if $A$ and $A^c$ are uniformly covered by $\set{\mathfrak{A}_n}_{n\in \nset{N}}$.

\end{theorem}

In view of the above theorem it is clear that it is convenient to establish the following definition. 

\begin{definition}\label{def:03}
Define $\mathfrak{A}_{\mu.a.e.}$ as 
$$
\mathfrak{A}_{\mu.a.e.}=\{A \in \mathfrak{A}: A \text{ and } A^c \text{ are uniformly covered by } \set{\mathfrak{A}_n}_{n\in \nset{N}} \}
$$
\end{definition}

Before we proceed we observe that if $\{A_n \in \mathfrak{A}_n\}$ is a sequence of sets that uniformly covers $A$ and  $\{C_n \in \mathfrak{A}_n\}$ is a sequence of sets that uniformly covers $A^c$, the sequence $A_n \bigcup C_n$ uniformly covers $A$ and $A^c$.

We have then that  $\mathfrak{A}_{\mu.a.e.}$ is an algebra of sets. 
XXXX No de hecho necesitamos que acomplemento este. Si ese es el caso si pero ¿por que?

Notice that if $\mathfrak{A}_{\mu.a.e.}$ is $\sigma$-subalgebra then

\begin{eqnarray*}
\underline{\mathfrak{A}}\subset \mathfrak{A}_{\mu.a.e.} \subset \mathcal{F}\subset \mathfrak{A}_\mu \subset \mathfrak{A}^\perp \subset \overline{\mathfrak{A}}.
\end{eqnarray*}

\begin{lemma}\label{lem:06}
If $\mathfrak{A}_{\mu.a.e.}$ is $\sigma$-subalgebra then

\begin{itemize}
\item[i)] $f\in L^\infty \left(\mathfrak{A}_{\mu.a.e.}\right)$ we have that

\begin{eqnarray*}
\mathcal{E}(f|\mathfrak{A}_n)\overset{a.e.}{\longrightarrow}\mathcal{E}(f|\mathcal{A}_{\mu.a.e.})=f.
\end{eqnarray*}

\item[ii)] If there is a $\sigma$-subalgebra $\mathfrak{B}$ such that for all $f\in L^\infty \left(\mathfrak{A}_{\mu.a.e.}\right)$

\begin{eqnarray*}
\mathcal{E}(f|\mathfrak{A}_n)\overset{a.e.}{\longrightarrow} \mathcal{E}(f|\mathfrak{B}) \text{ then}  \mathfrak{B} \subset\mathfrak{A}_{\mu.a.e.}.
\end{eqnarray*}

\end{itemize}

\begin{proof}

Let $\epsilon>0$. As $f\in L^\infty \left(\mathfrak{A}_{\mu.a.e.}\right)$, there is a simple function $g$, $\mathfrak{A}_{\mu.a.e.}$-measurable such that $\norm{f-g}_\infty<\epsilon/2$. It is clear that \textbf{Lemma \ref{lem:05}} implies that $\mathcal{E}(g|\mathfrak{A}_n)\overset{a.e.}{\longrightarrow} g$. Thus

\begin{align*}
\abs{\mathcal{E}(f|\mathfrak{A}_n)-f}(x)\leq & \abs{ \mathcal{E}(f|\mathfrak{A}_n)-\mathcal{E}(g|\mathfrak{A}_n) }(x)+\abs{\mathcal{E}(g|\mathfrak{A}_n)-g}(x)+\abs{g-f}(x)\\
\leq & \abs{\mathcal{E}(f-g|\mathfrak{A}_n) }(x)+\abs{\mathcal{E}(g|\mathfrak{A}_n)-g }(x)+\norm{g-f}_\infty \\
\leq& \epsilon +\abs{\mathcal{E}(g|\mathfrak{A}_n)-g}(x).
\end{align*}

So $\ds \varlimsup_{n\to \infty}\abs{\mathcal{E}(f|\mathfrak{A}_n)-f}(x)\leq \epsilon \hspace{0.1cm}a.e.$ for every $\epsilon$ and therefore $\ds \lim_{n\to \infty}\mathcal{E}(f|\mathfrak{A}_n)=f \hspace{0.1cm}a.e.$

To prove $ii)$, let $A\in \mathfrak{B}$.  By hypothesis   $\mathcal{E}(\chi_{A}|\mathfrak{A}_n)\longrightarrow \chi_A \hspace{0.1cm}a.e.$.

For $0<\epsilon<1$ define $A_n=\set{x \ : \ \mathcal{E}(\chi_A|\mathfrak{A}_n)>\epsilon } \in \mathfrak{A_n}$.

Since for almost all $x\in A$,  $\mathcal{E}(\chi_A|\mathfrak{A}_n)(x)\underset{n\to\infty}{\longrightarrow} 1$, $x\in A_n$ for $n$ big enough. So $\chi_{A}(x)\chi_{A_n}(x)\longrightarrow \chi_A(x)$ a.e..
The case for $A^c$ is similar. In this case, for almost  all $x\notin A$, $\mathcal{E}(\chi_A|\mathfrak{A}_n)(x)\underset{n\to\infty}{\longrightarrow} 0$. Hence  $x\notin A_n$ for $n$ big enough. Thus $\chi_{A^c}(x)\chi_{A_n}(x)\longrightarrow 0$ a.e.. And thus $\chi_{A_n}(x)\longrightarrow \chi_A(x)$ a.e..

\begin{eqnarray*}
A=\bigcup_{N=1}^\infty \bigcap_{n>N}A_n=\bigcap_{N=1}^\infty \bigcup_{n>N}A_n.
\end{eqnarray*}

we also have

\begin{align*}
\norm{A \setminus A_n}_{\mathfrak{A}_n}=\norm{\mathcal{E}(\chi_A \chi_{A_n^c}|\mathfrak{A}_n)}_\infty = \norm{\chi_{A_n^c}\mathcal{E}(\chi_A|\mathfrak{A}_n)  }_\infty < \epsilon \norm{\chi_{A_n^c}}_\infty <\epsilon.
\end{align*}

And so by \textbf{Lemma \ref{lem:09}} $A$ is uniformly covered by by $\set{\mathfrak{A}_n}$.

\end{proof}

\end{lemma}

We can improve the above lemma by considering functions in $L^p$.

\begin{lemma}\label{lem:06}
If $\mathfrak{A}_{\mu.a.e.}$ is $\sigma$-subalgebra and $1 \leq p < \infty$ then for all $f\in L^2 \left(\mathfrak{A}_{\mu.a.e.}\right)$ we have that

\begin{eqnarray*}
\mathcal{E}(f|\mathfrak{A}_n)\overset{a.e.}{\longrightarrow}\mathcal{E}(f|\mathcal{A}_{\mu.a.e.})=f.
\end{eqnarray*}

\end{lemma}

How do the above results look in the case of $\mathfrak{A}_n$ are monotone?. In the case we have a monotone decreasing sequence of $\sigma$-subalgebras  it is clear that $\mathfrak{A}_N=\bigvee_{n=N}^\infty \mathfrak{A}_n$ and that for any $N$,  $ \bigcap_{n=N}^\infty \mathfrak{A}_n =  \bigcap_{n=1}^\infty \mathfrak{A}_n.$ Therefore
$$
\bigcap_{n=1}^\infty \mathfrak{A}_n=\underline{\mathfrak{A}} \subset \mathcal{F} \subset \overline{\mathfrak{A}}\subset \bigcap_{n=1}^\infty \mathfrak{A}_n
$$

So if $A \in \mathcal(F) $ the sequence $A_n=A \in \mathfrak{A}_n $ and trivially uniformly covers $A$.

XXXXXX falta la demostracion

XXXXXXXXXXXXXXXXXXXXX

\subsection{Necessary and sufficient conditions to have}

In [4.x] we defined a $\sigma$-subalgebra $\mathfrak{A}^{\perp}$ by first considering the set $W$

\begin{eqnarray*}
W=\set{g\in L^2(\mathfrak{A}) \ : \ \exists A_{n_k}\in \mathfrak{A}_{n_k} \ \mbox{ with } \ \chi_{A_{n_k}} \overset{L^2-weakly}{\longrightarrow}g },
\end{eqnarray*}

 $\mathfrak{A}_\perp$ was defined as the minimum   $\sigma$-algebra  generated by $W$. We are going to do something similar. First we notice the obvious fact $\chi_{A_n}= \mathcal{E}(\chi_{A_n}|\mathfrak{A}_n)$. Define then,  
 given a sequence $\set{\mathfrak{A}_n}_{n\in \nset{N}}$ $\sigma$-subalgebras, the set

\begin{eqnarray*}
C_N=\set{h\in L^1(\mu) \ : h=\sum_{k\geq N}\mathcal{E}(\chi_{B_k}|\mathfrak{A}_k), \  B_k\in \mathfrak{A}, \text{ disjoint,  and such that there are only finite many of such }B_{k}'s  }.
\end{eqnarray*}

Notice that if $h\in C_N$

\begin{eqnarray*}
\norm{h}_1=\int \sum_{n\geq N}\mathcal{E}(\chi_{B_n}|\mathfrak{A}_n)d\mu= \sum_{n\geq N}\mu(B_n)= \mu\left( \bigcup_{n\geq N}B_n \right)\leq 1.
\end{eqnarray*}

XXX"""XXX"""

\
We will define
\begin{definition}

\begin{eqnarray*}
W^{\perp a.e.}=\set{f\in L^\infty(\mu)  \ : \ \mbox{for every subsequence }\set{h_{N_k}}, \ h_{N_k}\in C_{N_k}, \ \scp{f,h_{N_k}} \underset{k\to \infty}{\longrightarrow}0  }.
\end{eqnarray*}

\end{definition}

\begin{theorem}
If $f\in L^\infty(\mu)$ then

\begin{eqnarray}
\mathcal{E}(f|\mathfrak{A}_n)\overset{a.e.}{\longrightarrow}0,
\end{eqnarray}

if and only if $f\in W^{\perp a.e.}$.

\begin{proof}
$\Rightarrow)$
Let $f\in L^\infty(\mu)$. Without loss of generality we can assume that $\norm{f}_\infty\leq 1$. First, notice that if $h\in C_N$ we have

\begin{eqnarray*}
\scp{h,f}=\sum_{n\geq N}\scp{\mathcal{E}(\chi_{B_n},\mathfrak{A}_n),f}=\sum_{n\geq N}\scp{\chi_{B_n},\mathcal{E}(f|\mathfrak{A}_n)}.
\end{eqnarray*}

Let $\epsilon>0$. Since $\mathcal{E}(f,\mathfrak{A}_n)\overset{a.e.}{\longrightarrow}0$, Egoroff's theorem implies that there is a set $M_\epsilon$ and $N_1\in \nset{N}$ such that $m\geq N_1$

\begin{eqnarray*}
\mu(M_\epsilon^c)<\dfrac{\epsilon}{2} & \mbox{ and }&\norm{\chi_{M_\epsilon}\mathcal{E}(f|\mathfrak{A}_n)}<\dfrac{\epsilon}{2}.
\end{eqnarray*}

Thus if $N>N_1$

\begin{align*}
\abs{\scp{h_N,f}}&=\abs{\sum_{n\geq N}\scp{\chi_{B_n},\mathcal{E}(f,\mathfrak{A}_n)  }}\\
&\leq \sum_{n\geq N}\left(\abs{\scp{\chi_{B_n},\chi_{M_\epsilon}\mathcal{E}(f,\mathfrak{A}_n)  }}+ \abs{\scp{\chi_{B_n}\chi_{M_\epsilon^c},\mathcal{E}(f|\mathfrak{A}_n)  }} \right)\\
&\leq \dfrac{\epsilon}{2}\sum_{n\geq N}\scp{\chi_{B_n},\chi_{M_\epsilon}}+\sum_{n\geq N}\norm{\mathcal{E}(f|\mathfrak{A}_n)}_\infty \norm{\chi_{B_n}\chi_{M_\epsilon^c}}_1\\
&\leq \dfrac{\epsilon}{2}\mu(M_\epsilon)+\norm{f}_\infty\mu(M_\epsilon^c)\leq \epsilon,
\end{align*}

thus $\scp{h_N,f}\underset{N\to \infty}{\longrightarrow}0$.

$\Leftarrow)$ To prove the theorem in the other direction first we notice that since we can take $f$ or $-f$, without loss of generality we can assume that there is an $\epsilon$ such that

\begin{eqnarray*}
\mu\left(\bigcap_{N=1}^\infty \bigcup_{n\geq N} \set{x\in \nset{X} \ : \ \mathcal{E}(f|\mathfrak{A}_n)(x)\geq \epsilon} \right)>0.
\end{eqnarray*}

That is, there is an $r>0$ such that for any $N$ there is an $M$ with

\begin{eqnarray*}
\mu\left(\bigcup_{N\leq n\geq M} \set{x\in \nset{X} \ : \ \mathcal{E}(f|\mathfrak{A}_n)(x)\geq \epsilon} \right)>r.
\end{eqnarray*}

Let $A_n=\set{x\in \nset{X} \ : \ \abs{\mathcal{E}(f|\mathfrak{A}_n)(x)}>\epsilon}$ and as usual

\begin{eqnarray}
B_N=A_N, & B_k=A_k \setminus \bigcup_{j=N}^{k-1}B_j \ \mbox{ for }\ N<k<M,
\end{eqnarray}

$\set{B_k}$ is a disjoint family and

\begin{eqnarray*}
\bigcup_{N\leq n<M} B_n=\bigcup_{N\leq n<M}A_n.
\end{eqnarray*}

Let

\begin{eqnarray*}
h_N=\sum_{N\leq n<M}\mathcal{E}(\chi_{B_n}|\mathfrak{A}_n),
\end{eqnarray*}

we have that

\begin{eqnarray*}
\scp{h_N,f}=\sum_{N\leq n<M}\scp{\chi_{B_n}|\mathcal{E}(f|\mathfrak{A}_n)}>\epsilon \sum_{N\leq n<M}\mu(\chi_{B_n})=\epsilon \mu\left(\bigcup_{N\leq n<M} B_n\right)>\epsilon r.
\end{eqnarray*}

We have then constructed a sequence $\set{h_N}_{N\in \nset{N}}$ such that for all $N$,  $\scp{h_N,f}>\epsilon r$. Therefore $f\not\in W^{\perp a.e.}$.

\end{proof}

\end{theorem}

In XXX we defined the orthogonal conditional expectation induced by $\mathfrak{D}$ as the operator $\mathcal{E}_{\mathfrak{B}}^{\perp}=I-\mathcal{E}_{\mathfrak{B}}$, which in $L^2(\mathfrak{A})$ is the orthogonal projection $\mathcal{E}_{\mathfrak{B}}^{\perp}: L^2(\mathfrak{A}) \rightarrow L^2(\mathfrak{B})^\perp$. 
Let $\nset{D}$ be the following family of $\sigma$-subalgebras 
$$
\nset{D}=\set{\salg{B} \in \salg{A}: \cexpperp{\salg{B}}f \in W^{\perp a.e.} \text{ for all } f \in L^\infty (\salg{A})}
$$

It is clear that $\nset{D}$ is not empty since it trivially contains $\salg{A}$.

An immediate property is the following.
\begin{proposition}
Let $\salg{C}, \salg{B}$ be two $\sigma$-subalgebras. If $\salg{C} \in \nset{D}$ and $\salg{B} \supset \salg{C}$, then $\salg{B} \in \nset{D}$
\begin{proof}
Notice that in this case  $\cexpperp{\salg{B}}=\cexpperp{\salg{C}}\cexpperp{\salg{B}}$. So, as
 $f$ in $L^\infty$ implies that $\cexpperp{\salg{B}}f$ is also in $L^\infty$,
$$
\cexpperp{\salg{B}}f=\cexpperp{\salg{C}}(\cexpperp{\salg{B}}f) \in W^{\perp a.e.}
$$

\end{proof}
\end{proposition}
\begin{definition}
$\salg{A}_{min}$ will be the minimum complete $\sigma$-subalgebra such that contains the set $$\tilde{W}=\set{g \in L^2(\salg{A}): \text{there is } \set{h_{N_k}}, h_{N_k} \in C_{N_k}  \text{ such that }  h_{N_k} \longrightarrow g \text{ weakly in } L^2}$$ 

\end{definition}
\begin{lemma}
If 
$\salg{B} \in \nset{D}$ then $\salg{A}_{min} \subset \salg{B}$.
\end{lemma}

\begin{proof}
Let $g \in \tilde{W}$ and  $h_N \in C_N$ such that $ h_N \overset{w}{\longrightarrow} g$. Since $\salg{B}$ is in $\nset{D} $ we have that for all $f$ in $L^p \cap L^\infty$
$$
\scp{f,\cexp{g}{B}} = \scp{\cexp{f}{B},g} = \lim_{N\to \infty}\scp{\cexp{f}{B}, h_N}=\lim_{N\to \infty}\scp{(I-\cexpperp{\salg{B}}) f,h_N}=\lim_{N\to \infty}\scp{ f,(I-\cexpperp{\salg{B}}) h_N}
 = \lim_{N\to \infty}\scp{ f,h_N}=\scp{ f,g}
$$
Therefore $g$ is equal almost everywhere to $\cexp{g}{B}$ and so is $\salg{B}$ measurable. Since by definition $\salg{A}_{min}$ is the minumum $\sigma$-subalgebra such that makes all such $g$ measurable so  $\salg{A}_{min} \subset \salg{B}$.
\end{proof}
\begin{definition}
We will call the sequence $\set{\salg{A}_n}_{n\in \nset{N}} $ 2-bounded if for all $h_N$ in $C_N$ $\sup \norm{ h_N}_2 < \infty.$
\end{definition}

\begin{lemma}
If $\set{\salg{A}_n}$ is 2-bounded and $\salg{B} $ is a $\sigma$-subalgebra, then
$\salg{B} \in \nset{D}$ if and only if $\salg{A}_{min} \subset \salg{B}$.
\end{lemma}
\begin{proof}
In view of lemma xxx we only need to prove the assertion in the sense
$(\Leftarrow)$
Assume $\salg{A}_{min} \subset \salg{B}$ and suppose that $\salg{A}_{min}$ is not in $\nset{D}$. That means  that there is an $f \in L^\infty$ such that its orthogonal projection with respect to  $\salg{A}_{min}$ is not in $W^{\perp a.e.}$.That is, there is a sequence  
$\set{h_N \in C_N}$ such that $\scp{f,h_N}$ does not converge to zero. Hence, there is an $\epsilon >0$ and a subsequence $\set{h_{N_k}}$
such that $\abs{\scp{\cexpperp{\salg{A}_{min}}f,h_{N_k}}} > \epsilon$. By hypothesis the sequence of $\sigma$-subalgebras is 2-bounded, so there is a subsequence, which we still denote by $h_{N_k}$, that weakly converges to  an $h$ in $L^2$. By definition $h$ is in $\tilde{W}$ and therefore is $\salg{A_{min}}$ measurable. But that leads us to a contradiction since   $0 < \epsilon < \abs{\scp{\cexpperp{\salg{A}_{min}}f,h}} = \abs{\scp{f,\cexpperp{\salg{A}_{min}}h}} = 0$.

Of course the 2-boundeness condition is a very strong restriction. In the following we will show that the main property is the fact that $\nset{D}$ has a minimum $\sigma$-subalgebra.
In those cases we will denote that minimum by $\salg{A}_{\perp a.e.}$.

Of course in the case the sequence is 2-bounded 
$\salg{A}_{\perp a.e.} = \salg{A}_{min}$

\end{proof}
$$
$$
HOOOOOOOOOOOOOOOLA
XXXXXXXXXXXXXXXXXXXXX

\subsection{Necessary and sufficient conditions to have comvergence a.e.}


\begin{theorem}
If $\set{\salg{A}_n}$ is a sequence of $\sigma$-subalgebras such that

\begin{itemize}
\item[i)]
\item[ii)] 
\end{itemize}
for all $f \in L^{\infty}$, $\cexp{f}{A_n}$ converges a.e. then $\salg{A}_{\mu a.e.}$ is a $\sigma$-subalgebra, $\nset{D}$ has a minimum $\sigma$-subalgebra and $\salg{A}_{\mu a.e.} = \salg{A}_{\perp a.e.}$ 
\end{theorem}

\begin{proof}
By definition 
$\salg{A}_{\mu a.e.} = \{A \in \salg{A} | A \text{ and } A^c \text{ are univormly covered} \} $ Let us recall that for xxx $\salg{A}_{\mu a.e.}$ is an algebra. Thus, in order to prove that is a $\sigma-$subalgebra we need only to check that if we have a sequence of pairwise  disjoint sets $\{ B_k \}$, with
$ B_k \in \salg{A}_{\mu a.e.}$ , then  $\cup_{k=1}^\infty B_k \in \salg{A}_{\mu a.e.}$. 

In view of theorem 4.11, we need to prove that $\mathcal{E}(\chi_{\cup_{k=1}^\infty B_k}|\mathfrak{A}_n)\overset{a.e.}{\longrightarrow} \chi_{\cup_{k=1}^\infty B_k}. $ Now, since the $B_k$ are in $\salg{A}_{\mu a.e.}$ for any natural M we have
$$
\mathcal{E}(\chi_{\cup_{k=1}^\infty B_k}|\mathfrak{A}_n) \geq \mathcal{E}(\chi_{\cup_{k=1}^M B_k}|\mathfrak{A}_n)\overset{a.e.}{\longrightarrow} \chi_{\cup_{k=1}^M B_k}
$$
By hypothesis, the conditional expectations converge almost everywhere for any $f \in L^\infty$. Therefore there is a $g \in L^\infty$ such that
$$
\mathcal{E}(\chi_{\cup_{k=1}^\infty B_k}|\mathfrak{A}_n)\overset{a.e.}{\longrightarrow} g
$$
and thus $ g \geq  \chi_{\cup_{k=1}^M B_k}$. Since this is true for any M, we conclude that $ g \geq  \chi_{\cup_{k=1}^\infty B_k}.$ 

Using the dominated convergence theorem we have that
$$0 \leq \scp{1,g -\chi_{\cup_{k=1}^\infty B_k}} = \lim \scp{1,\mathcal{E}(\chi_{\cup_{k=1}^\infty B_k}|\mathfrak{A}_n)-\chi_{\cup_{k=1}^\infty B_k}} = \mu (\cup_{k=1}^\infty B_k) -\mu (\cup_{k=1}^\infty B_k) = 0
$$
$$0 \leq \int g -\chi_{\cup_{k=1}^\infty B_k} \, d \mu   = \lim \int \mathcal{E}(\chi_{\cup_{k=1}^\infty B_k}|\mathfrak{A}_n)-\chi_{\cup_{k=1}^\infty B_k} \, d\mu = \mu (\cup_{k=1}^\infty B_k) -\mu (\cup_{k=1}^\infty B_k) = 0
$$
And hence $ g =  \chi_{\cup_{k=1}^\infty B_k}$  a.e..

\end{proof}

\end{document}